\documentclass[12pt]{article}
\usepackage{amsmath, amsthm, amsfonts}

\newtheorem{theorem}{Theorem}
\newtheorem{corollary}[theorem]{Corollary}
\newtheorem{lemma}[theorem]{Lemma}
\newtheorem{proposition}[theorem]{Proposition}

\title{On the Areas of\\Cyclic and Semicyclic Polygons}
\author{F. Miller Maley$^*$ \and David P. Robbins$^*$ \and Julie Roskies$^*$}
\date{July 2004}

\def\qtr{{\textstyle\frac{1}{4}}}
\def\half{{\textstyle\frac{1}{2}}}
\def\abs#1{\left|#1\right|}
\def\floor#1{\lfloor#1\rfloor}

\def\Res{\mathop{\rm Res}}
\def\disc{\mathop{\rm discr}\nolimits}
\def\squish#1{\hbox to 0pt{\hss$\scriptstyle#1$\hss}}
\def\pile#1#2{\genfrac{}{}{0pt}{1}{#1}{#2}}
\def\Aspace{\mathbb{A}}
\def\Cfield{\mathbb{C}}
\def\Pspace{\mathbb{P}}
\def\Qfield{\mathbb{Q}}
\def\th{^{\rm th}}
\let\eps=\epsilon
\def\divides{\mathrel{|}}
\let\tild=\widetilde
\begin{document}
\maketitle
\footnotetext{* Center for Communications Research, Princeton, NJ 08540, USA
(maley@idaccr.org or julie@idaccr.org).}

\begin{abstract}
We investigate the ``generalized Heron polynomial'' that relates the
squared area of an $n$-gon inscribed in a circle to the squares of its
side lengths.  For a $(2m+1)$-gon or $(2m+2)$-gon, we express it as the
defining polynomial of a certain variety derived from the variety of binary
$(2m-1)$-forms having $m-1$ double roots.  Thus we obtain explicit
formulas for the areas of cyclic heptagons and octagons, and illuminate
some mysterious features of Robbins' formulas for the areas of cyclic
pentagons and hexagons.  We also introduce a companion family of
polynomials that relate the squared area of an $n$-gon inscribed in a
circle, one of whose sides is a diameter, to the squared lengths of the
other sides.  By similar algebraic techniques we obtain explicit
formulas for these polynomials for all $n\le 7$.
\end{abstract}

\section{Introduction}\label{s:intro}

Heron of Alexandria (c. 60 BC) is credited with the formula that relates
the area $K$ of a triangle to its side lengths $a$, $b$, and $c$:
$$K = \sqrt{s(s-a)(s-b)(s-c)}$$ where $s=(a+b+c)/2$ is the
semiperimeter.  For polygons with more than three sides, the side
lengths do not in general determine the area, but they do if the polygon
is convex and \emph{cyclic} (inscribed in a circle).  Brahmagupta,
in the seventh century, gave the analogous formula for a convex cyclic quadrilateral
with side lengths $a$, $b$, $c$, and $d$:
$$K = \sqrt{(s-a)(s-b)(s-c)(s-d)}$$ where $s=(a+b+c+d)/2$.  See
\cite{coxeter-greitzer} for an elementary proof.

Robbins \cite{robbins} found a way to generalize these formulas.  First,
drop the requirement of convexity and consider the square of the
(signed) area $K$ of a possibly self-intersecting oriented cyclic
polygon.  For this purpose we can define the area enclosed by a closed
curve $C$ to be $\oint_C x\,dy$.  Second, express the relation between
$K^2$ and the side lengths as a polynomial equation with integer
coefficients.  Given a cyclic polygon, one can permute its edges within
its circumscribed circle without changing its area, so the polynomial
will be symmetric in the side lengths, and in fact it can be written in
terms of $16K^2$ and the elementary symmetric functions $\sigma_i$ in
the \emph{squares} of the side lengths.  For instance, the Heron and
Brahmagupta formulas can be written
$$16K^2 - 4\sigma_2 + \sigma_1^2 - \eps\cdot 8\sqrt{\sigma_4} = 0$$
in which $\eps$ is $0$ for a triangle, $1$ for a convex
quadrilateral, and $-1$ for a nonconvex quadrilateral.  Hence
all cyclic quadrilaterals satisfy the polynomial equation
$(16K^2 - 4\sigma_2 + \sigma_1^2)^2 - 64\sigma_4 = 0$.  The general
result is the following.

\begin{theorem}\label{t:cycpoly}\rm\cite{robbins}\ \sl
For each $n\ge 3$, there is a unique (up to sign) irreducible
polynomial $\alpha_n$ with integer
coefficients, homogeneous in $n+1$ variables with the first variable
having degree~2 and the rest having degree~1, such that $\alpha_n(16K^2, a_1^2,
\ldots, a_n^2)=0$ whenever $a_1$, \dots, $a_n$ are the side lengths
of a cyclic $n$-gon and $K$ is its area.
\end{theorem}

The polynomials $\alpha_n$ are now known in the literature as
\emph{generalized Heron polynomials}.  For certain sets of $n$ side
lengths, as shown in \cite{robbins}, one can find up to $\Delta_n$
distinct squared areas, where
$$\Delta_n = \frac{n}{2} \binom{n-1}{\floor{\frac{n-1}{2}}} - 2^{n-2}.$$
Hence one expects that $\alpha_n$ has degree $\Delta_n$ in its first
variable.  This conjecture of Robbins, and two others made in
\cite{robbins}, have recently been established.  We summarize them in
Theorem~\ref{t:russians}.

\begin{theorem}\label{t:russians}\sl
The polynomial $\alpha_n$ is monic in $16K^2$ and has total degree
$2\Delta_n$.  If $n$ is even, then $\alpha_n = \beta_n \beta_n^*$ where
$\beta_n$ is a polynomial in $16K^2$, $\sigma_1$, \dots, $\sigma_{n-1}$,
$\sqrt{\sigma_n}$, where $\sqrt{\sigma_n} = a_1\cdots a_n$ and $\beta_n^*$
is $\beta_n$ with $\sqrt{\sigma_n}$ replaced by $-\sqrt{\sigma_n}$.
\end{theorem}

See \cite{fedorchuk-pak} or \S\ref{s:radius} for the degree,
\cite{varfolomeev} or \cite{connelly} for monicity, and
\cite{varfolomeev} for the factorization when $n$ is even.  Robbins'
main interest, however, and the motivation for our research, was to find
reasonably explicit formulas for all~$\alpha_n$ and $\beta_n$.


In \cite{robbins}, Robbins found formulas for $\alpha_5$ and $\beta_6$
that have a curious form.  To present them concisely, we introduce the
\emph{crossing parity} $\eps$ of a cyclic $n$-gon.  Assume that the
$n$-gon has vertices $v_1$, \dots, $v_n$ in the complex plane and
circumcenter $0$.  For odd $n$ let $\eps=0$, and for even $n$ let $\eps
= \pm 1$ be
\begin{equation*}
\eps = \mathop{\rm sign} \bigl[(1-v_2/v_1)(1-v_3/v_2)\cdots(1-v_1/v_n)\bigr],
\end{equation*}
which is well-defined because the product is real and nonzero if all
edges have positive length.  (This definition agrees with the previous
definition of $\eps$ for $n\in\{3,4\}$.)  Now assume $n\in \{5,6\}$.
Define $u_2 = -4K^2$, and make the substitutions
\begin{equation}\label{e:subs}
\begin{aligned}
t_1 &= \phantom{-}\sigma_1,\\
t_2 &= -\sigma_2 + \qtr t_1^2\phantom{t_2} - u_2,\\
t_3 &= \phantom{-}\sigma_3 + \half t_1 t_2 - \eps\cdot 2\sqrt{\sigma_6},\\
t_4 &= -\sigma_4 + \qtr t_2^2\phantom{t_1} + \eps\cdot t_1\sqrt{\sigma_6},\\
t_5 &= \phantom{-}\sigma_5 \mathbin{\phantom{+}} \phantom{\half t_1 t_2}
+ \eps\cdot t_2\sqrt{\sigma_6}.
\end{aligned}
\end{equation}
Then, for any cyclic pentagon or hexagon of the given crossing parity,
the cubic polynomial $u_2 + t_3 z + t_4 z^2 + t_5 z^3$ has a double
root, so its discriminant vanishes:
\begin{equation*}
t_3^2 t_4^2 - 4 u_2 t_4^3 - 4 t_3^3 t_5 + 18 u_2 t_3 t_4 t_5 - 27
u_2^2 t_5^2 = 0.
\end{equation*}
When the $t_i$ are expanded, this discriminant is a polynomial of degree
$\Delta_5 = 7$ in $u_2$, and hence in $16K^2$.  Multiplying it by
$2^{18}$ makes it monic in $16K^2$ and yields $\alpha_5$, $\beta_6$, or
$\beta^*_6$ according to whether $\eps$ is $0$, $+1$, or~$-1$.

In \S\ref{s:general}, we generalize this construction.  Fix $n$ and
the crossing parity~$\eps$, and let $m=\floor{(n-1)/2}$.  We introduce
auxiliary quantities $u_2$, \ldots, $u_m$, with $u_2 = -4K^2$, and
inductively define certain polynomial expressions $t_i$ in the
$\sigma_j$ and $u_j$ with~$j\le i$.  For $n=5$ or $6$, these definitions
reduce to (\ref{e:subs}).  Corollary~\ref{c:cyclic} then says that
the polynomial
$$P_n(z) = u_2 + \cdots + u_m z^{m-2} + t_{m+1} z^{m-1} + \cdots +
t_{2m+1} z^{2m-1}$$ is divisible by the square of a
polynomial of degree~$m-1$.  In other words, for any values of the $t_i$
and $u_j$ coming from a cyclic $n$-gon, $P_n(z)$ has $m-1$ double roots
over $\Cfield$
(counting with multiplicity, and including roots at infinity).  Such
polynomials form a variety of codimension $m-1$, defined locally by
$m-1$ equations.  So, if we regard $u_3$ through $u_m$ as indeterminates
and expand each $t_i$ in terms of the $\sigma_j$ and $u_j$, we can
in principle
eliminate the $m-2$ unwanted quantities $u_3$, \dots,~$u_m$ and recover
a single polynomial, which is $\alpha_{2m+1}$, $\beta_{2m+2}$, or
$\beta^*_{2m+2}$ depending on~$\eps$.  In \S\ref{s:formulas} we carry
out this program for $m=3$ to obtain formulas for $\alpha_7$ and
$\alpha_8$, the generalized Heron polynomials for cyclic heptagons and
octagons.

There is another family of area polynomials susceptible to the same
analysis.  Call an $(n+1)$-gon \emph{semicyclic} if it is inscribed in a
circle with one of its sides being a diameter.  Its squared area
satisfies a polynomial relation with the squares of the lengths of the
other $n$ sides; the degree in the squared area turns out to be
$$\Delta'_n = \frac{n}{2} \binom{n-1}{\floor{\frac{n-1}{2}}} =
\Delta_n + 2^{n-2}.$$

\begin{theorem}\label{t:semipoly}\sl
For each $n\ge 2$, there exists a unique monic irreducible
polynomial $\alpha'_n$ with integer coefficients, homogeneous in
$n+1$ variables with the first variable having
degree~2 and the rest having degree~1, such that $\alpha'_n(16K^2, a_1^2,
\ldots, a_n^2) = 0$ whenever $a_1$, \dots, $a_n$ are the lengths of the
sides of a semicyclic $(n+1)$-gon excluding a diameter, and $K$ is its
area.  The total degree of $\alpha_n'$ is $2\Delta'_n$.
\end{theorem}

The proof that $\alpha'_n$ exists and is unique (without assuming
monicity) follows the proof of Theorem~\ref{t:cycpoly} in
\cite{robbins} almost verbatim, and the argument in \cite{connelly}
shows that $\alpha'_n$ is monic.  We establish the degree by an
elementary argument in \S\ref{s:radius}, which is independent of the
rest of this paper.

Cyclic and semicyclic polygons are similar in many ways.  For instance,
just as the polygon of largest area one can make with $n$ given side
lengths is convex and cyclic, the polygon of largest area one can make
with $n$ given side lengths and one free side is convex and semicyclic.
We will adduce many algebraic similarities in the following sections.
For now we just observe that the polynomial $\alpha'_3$, which can be
worked out by hand, also takes the form of a discriminant: if
$u_2 = -4K^2$, then
\begin{equation}\label{e:semiquad}
\alpha'_3 = 16 \disc_z \bigl( z^3 + \sigma_1 z^2 + (\sigma_2 + u_2) z +
\sigma_3 \bigr).
\end{equation}


\section{The Main Identity}\label{s:main}

All our area formulas are based on a generating function identity that
relates the symmetric functions $\sigma_i$ in the squared side lengths
to certain quantities $\tau_j$ that arise in Robbins' proofs of the
pentagon and hexagon formulas.  The identity, Theorem~\ref{t:main},
holds for both cyclic and semicyclic polygons and for both odd and
even~$n$.

Suppose we have a cyclic $n$-gon or semicyclic $(n+1)$-gon inscribed in
a circle of radius~$r$ centered at the origin in the complex plane.  Let
its vertices be $v_1$, \dots, $v_n$ and $v_{n+1} = \delta v_1$, where
$\delta = 1$ for a cyclic $n$-gon and $\delta=-1$ for a semicyclic
$(n+1)$-gon.  Introduce the vertex quotients $q_i = v_{i+1}/v_i$ for
$i=1$, \dots,~$n$, and let $\tau_0$, $\tau_1$, \dots,~$\tau_n$ be the
elementary symmetric functions of the~$q_i$.  Then $\tau_0 = 1$ and
$\tau_n = q_1 q_2 \cdots q_n = \delta$.  Elementary geometry yields
the equations
\begin{align}
\label{e:vq} a_i^2 &= r^2(2 - q_i - q_i^{-1}),\qquad 1\le i\le n,\\
\label{e:area} 16K^2 &= -r^4(q_1 + \cdots + q_n - q_1^{-1} - \cdots
- q_n^{-1})^2\\ &\notag = -r^4(\tau_1 - \delta \tau_{n-1})^2.
\end{align}
Using (\ref{e:vq}) one can express each $\sigma_i$ in terms of $r$ and
the~$\tau_i$.  Let $g(y) = y^2 + (x/r^2 - 2) y + 1$.  Observe that $x$
is one of the values $a_i^2$ exactly when $g(y)$ has one of the vertex
quotients $q_i$ as a root, or in other words, when $g(y)$ has a common
root with the polynomial $f(y) = \prod_{i=1}^n (y - q_i) = \sum_{i=0}^n
(-1)^i \tau_i y^{n-i}.$ Hence the resultant of $f(y)$ and $g(y)$ is a
constant times $$h(x) = \prod_{i=1}^n (x - a_i^2) = \sum_{i=0}^n (-1)^i
\sigma_i x^{n-i},$$ and the coefficient of $x^n$ reveals that the
constant is $\delta r^{-2n}$.  By expanding the resultant, one finds
that each $\sigma_i$ is $r^{2i}$ times a quadratic polynomial in
the~$\tau_i$.  A particularly simple example is
\begin{equation*}
\sigma_n = \delta (-1)^n r^{2n} (\tau_0 - \tau_1 + \tau_2 - \cdots \pm
\tau_n)^2.
\end{equation*}
If $n$ is even and $\delta=1$, then $\sqrt{\sigma_n}$ is expressible in
terms of $r^2$, the~$\tau_i$, and the crossing parity~$\eps$:%
\begin{equation}\label{e:sqrt}
\begin{aligned}\sqrt{\sigma_n} &= \abs{v_1 - v_2} \cdots \abs{v_n - v_{n+1}}
= r^n \abs{1 - q_1} \cdots \abs{1 - q_n}\\
&= r^n \eps (1 - q_1) \cdots (1 - q_n)
= r^n \eps (\tau_0 - \tau_1 + \tau_2 - \cdots + \tau_n).
\end{aligned}
\end{equation}
So far we are following \cite{robbins} except for the addition of the
semicyclic case.

Consider now the involution that reflects the polygon in the real axis.
This operation preserves the squared area and the side lengths, but it
replaces each $q_i$ with $q_i^{-1}$ and hence replaces each $\tau_i$
with $\delta\tau_{n-i}$.  Because each $\sigma_i$ is a quadratic form in
the $\tau_j$ preserved by the involution, it can be uniquely decomposed
into a two parts: a quadratic form in symmetric linear combinations of
the $\tau_j$, and a quadratic form in antisymmetric linear combinations
of the~$\tau_j$.  When we perform this decomposition on the whole
generating function $\sum_i (-x)^i \sigma_i$, each part factors in a
surprising way, which our main identity records.

To write the identity explicitly, we need the following linear
combinations of the $\tau_j$, for $0\le k\le n/2$:
\begin{align}\label{e:defd}
d_k &= \sum_{i=0}^k (-1)^i \binom{n-1-2k+i}{i} (\tau_{k-i} - \tau_{n-k+i}),\\
\label{e:defe}
e_k &\span {}= \sum_{i=0}^k (-1)^i \left[ \binom{n-2k+i}{i} +
\binom{n-2k+i-1}{i-1} \right] (\tau_{k-i} + \tau_{n-k+i}).
\end{align}
Let $D(x) = \sum d_i x^i$ and $E(x) = \sum e_i x^i$.

\begin{theorem}[Main Identity]\label{t:main}\sl
For a cyclic $n$-gon or semicyclic $(n+1)$-gon of radius~$r$, with
$\delta=1$ or $-1$ respectively, the symmetric functions $\sigma_i$ of
the squared side lengths and $\tau_i$ of the vertex quotients are
related by
$$\delta \cdot \sum_{i=0}^n (-x)^i \sigma_i = \qtr E(r^2x)^2 + (r^2 x -
\qtr) D(r^2 x)^2.$$
\end{theorem}

\begin{proof} When the $\sigma_i$ are expanded in terms of the $\tau_j$,
both sides of the main identity become polynomials in $r^2 x$, so we may
assume $r=1$.  The left-hand side is then 
$$\delta \cdot \sum_{i=0}^n (-x)^i \sigma_i = \delta x^n h(x^{-1}) = x^n \Res
( f, g )$$ where $f(y) = \sum_{i=0}^n (-1)^i \tau_i y^{n-i}$ and
$g(y) = (y-1)^2 + x^{-1}y$.

We calculate the resultant using its $PGL(2)$-invariance and other
standard properties \cite{g-k-z}.  Make the change of variable
$y = (z-1)/(z+1)$ so that the roots of $g$ are related
by $z \mapsto -z$ instead of $y \mapsto y^{-1}$.  We obtain the
polynomials
\begin{align*}
f^*(z) &= (z+1)^n f \left(\textstyle \frac{z-1}{z+1}\right) =
\sum_{i=0}^n (-1)^i \tau_i (z-1)^{n-i} (z+1)^i,\\
g^*(z) &= (z+1)^2 g \left(\textstyle \frac{z-1}{z+1}\right) =
x^{-1} (z^2 + 4x-1),
\end{align*}
and the transformation has determinant $2$, so
$$x^n \Res(f,g) = 2^{-2n} x^n \Res(f^*,g^*) = 2^{-2n} f^* \left( \!
\sqrt{1-4x} \right) f^* \left(\!\ -\sqrt{1-4x} \right).$$
Write $f^*(z) = f_0 (z^2) + z f_1(z^2)$, separating even and odd powers of~$z$.
Then $$x^n \Res(f, g) = 2^{-2n} \left[ \textstyle f_0 \bigl( 1-4x \bigr)^2 -
(1-4x) f_1 \bigl( 1-4x \bigr)^2 \right],$$ which explains the form of the main
identity.  It remains to evaluate $f_0(1-4x)$ and $f_1(1-4x)$.
We consider only $f_0$, as $f_1$ is similar but simpler.

It helps to introduce the Fibonacci polynomials
$F_n(x) = \sum_{i=0}^{\floor{n/2}} \binom{n-i}{i} x^i$,
which count compositions of $n$ by $1$'s and $2$'s.
They satisfy the recurrence
\begin{equation}\label{e:recur}
F_n(x) = F_{n-1}(x) + x F_{n-2}(x), \qquad n\ge 1,
\end{equation}
with $F_0(x) = 1$ and $F_n(x)=0$ for $n<0$, and have generating function
$$F(x;t) = \sum_n F_n(x) t^n = (1 - t - xt^2)^{-1}.$$
The generating functions for the $d_k$ and $e_k$ can then be written
\begin{align*}
D(x) &= \sum_{i=0}^{\floor{n/2}}
(\tau_i - \tau_{n-i}) x^i F_{n-2i-1}(-x),\\
E(x) &= \sum_{i=0}^{\floor{n/2}}
(\tau_i + \tau_{n-i}) x^i \bigl( F_{n-2i}(-x) - x F_{n-2i-2}(-x) \bigr).
\end{align*}

To evaluate $f_0(1-4x)$, first rewrite $f_0(z^2) = \half \bigl( f^*(z) + f^*(-z) \bigr)$ in
terms of the sums $\tau_i + \tau_{n-i}$.  If we let $\theta=\half$ if $2i=n$ and
$\theta=1$ otherwise, then
$$f_0(z^2) = \sum_{i=0}^{\floor{n/2}} (-1)^{n-i} \theta (\tau_i + \tau_{n-i})
(1-z^2)^i \sum_{j} \binom{n-2i}{2j} z^{2j}.$$
Next, we want to substitute $z^2 = 1-4x$ and evaluate the sum over $j$.
Writing $m=n-2i$, we can compute the generating function
$$\sum_{m\ge 0} t^m \sum_j \binom{m}{2j} (1-4x)^j = \frac{1-t}{1 - 2t + 4xt^2}
= (1-t) F(-x;2t)$$ by interchanging sums and simplifying.  Thus we get
$$f_0(1-4x) = \sum_{i=0}^{\floor{n/2}} (-1)^{i+n} \theta (\tau_i + \tau_{n-i})
(4x)^i 2^{n-2i} \left( F_{n-2i}(-x) - \half F_{n-2i-1}(-x) \right).$$
The recurrence (\ref{e:recur}) shows that $\theta \left(
F_{n-2i}(-x) - \half F_{n-2i-1}(-x) \right)$ is equivalent to
$\half \left(F_{n-2i}(-x) - x F_{n-2i-2}(-x) \right)$, and so
$$f_0(1-4x) = (-1)^n 2^{n-1} E(x).$$
Likewise $f_1(1-4x) = (-1)^n 2^{n-1} D(x)$, and the identity follows.
\end{proof}

\section{Consequences of the Main Identity}\label{s:general}

The main identity tells us how to generalize the definition of the
quantities $t_i$ and $u_j$ that were so useful in simplifying the
pentagon and hexagon formulas.  Cyclic $n$-gons have $d_0 = \tau_0 -
\tau_n = 0$ and $e_0 = \tau_0 + \tau_n = 2$ by (\ref{e:defd}) and
(\ref{e:defe}), so the expansions of $E(r^2 x)^2$ and $D(r^2 x)^2$
include linear terms in the $e_k$ but not in the~$d_k$.  The
substitutions that replace the $\sigma_i$ with the $t_i$ will first
isolate and then eliminate the variables $e_k$, and leave us with
expressions relating the $t_i$ and $u_j$ to the radius $r$ and the
variables $d_k$.  The main identity will then express the algebraic
relationship among the $t_i$ and $u_j$ as the factorization of a single
polynomial $P_n(z)$.

\begin{corollary}\label{c:cyclic}\sl
Given a cyclic $n$-gon of crossing parity $\eps$ and radius $r$, let
$m=\floor{(n-1)/2}$ and let $u_j = r^{2j} \sum_{i=1}^{j-1} (\qtr d_i -
d_{i-1}) d_{j-i}$ for $j\ge 1$.  Inductively define $t_0 = -2$ and
$$t_j = (-1)^{j+1} \sigma_j + \sum_{1\le i,j-i \le m}
\frac{t_i t_{j-i}}{4} + \begin{cases} -u_j,& \text{if $j\le m$,}\\
\eps \cdot (-1)^m t_{j-m-1} \sqrt{\sigma_n}, & \text{if $j>m$,}\end{cases}$$
for $j=1, \ldots, 2m+1$.  Then $t_j = -e_j r^{2j}$ for $0\le j\le m$,
and the polynomial
$$P_n(z) = u_2 + u_3 z + \cdots + u_m z^{m-2} + t_{m+1} z^{m-1} +
\cdots + t_{2m+1} z^{2m-1}$$ factors as $(\qtr-r^2 z) [z^{-1} D(r^2 z)]^2$.
\end{corollary}

By (\ref{e:defd}) and (\ref{e:area}), we have $u_2 = \qtr r^4 d_1^2 = \qtr r^4
(\tau_1 - \tau_{n-1})^2 = -4K^2$, and $u_1 = 0$ by definition.  Thus the
$t_j$ and $u_j$ in Corollary~\ref{c:cyclic} agree with those defined
in \S\ref{s:intro}.

\begin{proof}[Proof of Corollary \ref{c:cyclic}]
We have $e_0 = \tau_0 + \tau_n = 2$ by the definition (\ref{e:defe}),
so $t_0 = -e_0 r^0$.  Now let $1\le j\le m$.  The coefficient of $x^j$
in $\qtr E(r^2x)^2$ is
$$r^{2j} \sum_{i=0}^{j} \frac{e_i e_{j-i}}{4} = r^{2j} e_j +
\sum_{i=1}^{j-1} \frac{t_i t_{j-i}}{4}$$ by induction on $j$.
The coefficient of $x^j$ in $(r^2 x-\qtr) D(r^2 x)^2$ is $-u_j$, so the
equation $t_j = -e_j r^{2j}$ follows by comparing coefficients of $x^j$
in the main identity.

For $j>m$ we must consider the coefficient $r^{2m+2} e_{m+1}$ of
$x^{m+1}$ in $E(r^2 x)$.  If $n=2m+1$, then $E$ has degree $m$ by
definition so this coefficient is zero.  But if $n=2m+2$, then the
coefficient is
\begin{align*}
r^{2m+2} e_{m+1} &= r^n \biggl( 2\tau_{m+1} + \sum_{i=1}^{m+1} (-1)^i
2(\tau_{m+1-i}+\tau_{m+1+i}) \biggr)\\
&= 2(-1)^{m+1} \eps \sqrt{\sigma_n}
\end{align*}
by (\ref{e:defe}) and (\ref{e:sqrt}).  So for $m<j\le 2m+1$, the
coefficient of $x^j$ in $\qtr E(r^2 x)^2$ is
$$r^{2j} \sum_{i=j-m}^m \frac{e_i e_{j-i}}{4} + r^{2j} \frac{e_{j-m-1}
e_{m+1}}{2} = \sum_{i=j-m}^m \frac{t_i t_{j-i}}{4} + t_{j-m-1} (-1)^m
\eps \sqrt{\sigma_n},$$ and this holds whether $n$ is odd or even
because $\eps=0$ when $n$ is odd.  Thus, by the main identity, $-t_j$ is
the coefficient of $x^j$ in $(r^2 x-\qtr)D(r^2 x)^2$ for $j=m+1$, \dots,
$2m+1$.

We now see that $(r^2 x - \qtr) D(r^2 x)^2$, a polynomial of degree
$2m+1$ whose two lowest terms vanish, is exactly $- x^2 P_n(x)$.
\end{proof}

There is a geometric argument that Corollary~\ref{c:cyclic} contains
enough information to recover $\alpha_n$.  To simplify the explanation,
assume $n=2m+1$ and $m\ge 2$.  The nonzero polynomials of degree up to
$2m-1$ that have a squared factor of degree $m-1$ naturally form a
projective variety of codimension $m-1$ in $\Pspace^{2m-1}$, which is
irreducible because it is the image of $\Pspace^1 \times \Pspace^{m-1}$
under a regular map.  Hence the affine variety $X_m \subset \Aspace^{2m}$
of such polynomials (now including the zero polynomial), which has the
same ideal, is also irreducible.  The substitutions that write
$t_{m+1}$, \dots, $t_{2m+1}$ in terms of the $\sigma_i$ and $u_j$ amount
to a morphism $f:\Aspace^{3m}\to \Aspace^{2m}$, which is a product
bundle with fiber $\Aspace^m$.  This is because, for any point $(u_2$,
\dots, $u_m$, $t_{m+1}$, \dots, $t_{2m+1})$ in the range and any given
values of $\sigma_1$, \dots, $\sigma_m$, the values of $\sigma_{m+1}$,
\dots, $\sigma_{2m+1}$ are uniquely determined as polynomial functions
of the other variables.  Hence $f^{-1}(X_m) \approx X_m \times \Aspace^m$ is
irreducible and of codimension $m-1$.  Finally, when we apply the
projection $\pi:\Aspace^{3m} \to \Aspace^{2m+2}$ that eliminates the
$m-2$ variables $u_3$, \dots, $u_m$, the closure of the image
$\pi(f^{-1}(X_m))$ is an irreducible variety of codimension at least $1$
that contains $V(\alpha_n)$, so it must equal $V(\alpha_n)$.

Corollary~\ref{c:cyclic} therefore reduces the problem of finding $\alpha_n$
to two subproblems: finding the defining equations of the
variety~$X_m$, and then, after expanding $t_{m+1}$, \dots, $t_{2m+1}$ in
terms of the $\sigma_i$ and $u_j$, eliminating the $m-2$ variables
$u_3$, \dots,~$u_m$.

The application of the main identity to semicyclic polygons is similar
but slightly different.  In this case $e_0 = \tau_0 + \tau_n = 0$ and
$d_0 = \tau_0 - \tau_n = 2$, so the main identity involves linear terms
in the $d_k$ but not the $e_k$.  Our definitions of $t_i$ and $u_j$ are
therefore designed to extract and eliminate the variables $d_k$.  Again
we can distill the relationship among the $t_i$ and $u_j$ to the
factorization of a polynomial $P'_n(z)$.  This time, due to the factor
$(r^2 x - \qtr)$ in the main identity, the expression for $t_i$
explicitly includes~$r^2$, so there remains one more unwanted variable
to eliminate for a given~$n$.

\begin{corollary}\label{c:semi}\sl
Given a semicyclic $(n+1)$-gon of radius $r$, let $m=\floor{(n-1)/2}$
and let $u_j = r^{2j} \sum_{i=1}^{j-1} e_i e_{j-i}/4$ for $1\le j\le m$.
Inductively define $t_0 = -2$ and
\begin{equation*}
t_j = (-1)^{j+1} \sigma_j + \sum_{\pile{1\le i\le m}{1\le j-i \le m}}
\frac{t_i t_{j-i}}{4} - r^2 \sum_{\pile{0\le i-1\le m}{0\le j-i\le m}}
t_{i-1} t_{j-i} + \begin{cases} -u_j,& \text{if $j\le m$,}\\ 0,&
\text{if $j>m$,}\end{cases}
\end{equation*}
for $j=1, \ldots, n$.  Then $t_j = -d_j r^{2j}$ for $0\le j\le m$, and
the polynomial
$$P'_n(z) = u_2 + u_3 z + \cdots + u_m z^{m-2} + t_{m+1} z^{m-1}
+ \cdots + t_{n} z^{n-2}$$
is the square of $E(r^2 z)/2z$.  In particular, $t_n=0$ if $n$ is odd.
\end{corollary}

Once again $u_2 = \qtr r^4 e_1^2 = \qtr r^4 (\tau_1 + \tau_{n-1})^2 = -4K^2$ by
(\ref{e:area}), since now $\delta=-1$.

\begin{proof} As in Corollary~\ref{c:cyclic}, the claims follow from
equating coefficients of $x^j$ in the main identity and inducting on~$j$
to evaluate $t_j$ for $0\le j\le m$.  If $n=2m+1$, the degree of $E(x)$ is
just $m$, so $t_n = 0$.
\end{proof}

The polynomial $P'_n(z)$ contains $m-1$ unwanted variables, namely $r^2$
and $u_3$, \dots, $u_m$.  If $n=2m+1$, then $P'_n(z)$ is a polynomial of
degree $2m-2$ that is a square, which gives rise to $m-1$ equations in
its coefficients, and we have the additional equation $t_n = 0$.  If
$n=2m+2$, then $P'_n(z)$ is a square of degree $2m$, which yields $m$
equations.  In either case Corollary~\ref{c:semi} holds enough
information, in principle, to derive the area formula $\alpha'_n$.  As
before, one can make this claim precise using some algebraic geometry.

\section{Explicit Formulas}\label{s:formulas}

In this section we apply the results of \S\ref{s:general} to produce
area formulas for cyclic heptagons and octagons, and also semicyclic
quadrilaterals, pentagons, hexa\-gons, and heptagons.

Because the degree of the generalized Heron polynomial $\alpha_n$ is
exponential in $n$, and the number of terms could be exponential in
$n^2$, there is some question as to what constitutes an explicit
formula.  Our formulas have concise descriptions, and if a
polygon is given with exact (for instance, rational) side lengths, the
polynomial satisfied by its area can be computed exactly using standard
operations such as evaluating determinants.

First let us apply Corollary~\ref{c:cyclic} to the cases
$n\in\{7,8\}$.  It gives us a binary quintic form
$$x^5 P_n(y/x) = u_2 x^5 + u_3 x^4 y + t_4 x^3 y^2 + t_5 x^2 y^3 + t_6 x
y^4 + t_7 y^5$$ whose coefficients are polynomials in $u_2$, $u_3$,
$\sigma_1$, \dots, $\sigma_7$ and perhaps $\sqrt{\sigma_8}$, and which,
when its coefficients are evaluated for any cyclic $n$-gon, has two
linear factors over $\Cfield$ of multiplicity two.  The condition for a
quintic form $Q$ to factor in this way is given by the vanishing of a
certain covariant $C$, which in the notation of transvectants
\cite{olver} is
$$C = 2 Q (H,i)^{(2)} + 25 H (Q,i)^{(2)} + 6 Q i^2, \quad H =
(Q,Q)^{(2)}, \quad i = (Q,Q)^{(4)}.$$ Here $(f,g)^{(d)} = \sum_{i=0}^d
(-1)^i \binom{d}{i} (\partial^d f / \partial x^i \partial y^{d-i})
(\partial^d g / \partial x^{d-i} \partial y^i).$ This fact about
quintics is presumably classical, but we have not yet found a reference.
In any case, $C$ is a form of degree~$9$ in $\{x,y\}$ whose coefficients
are forms of degree~$5$ in the coefficients of the original quintic, so
its coefficients give us ten degree-5 polynomials in $u_2$, $u_3$,
$t_4$, $t_5$, $t_6$, $t_7$ that must vanish.  These same ten polynomials
can be obtained as the Gr\"obner basis, with a graded term ordering, for
the ideal of the variety of quintic forms that factor as a linear form
times the square of a quadratic.


To obtain the desired relation between $u_2$ and the $\sigma_i$, we must
expand the coefficients of $C$ as polynomials in $u_3$ and then
eliminate $u_3$.  We can do this most explicitly using resultants with
respect to $u_3$.  The two simplest coefficients of $C$ are
\begin{align*}
F &= u_3^2 t_4^3 - 4 u_2 t_4^4 - 4 u_3^3 t_4 t_5 + 18 u_2 u_3 t_4^2 t_5
- 27 u_2^2 t_4 t_5^2\\ &\qquad + (8 u_3^4 - 42 u_2 u_3^2 t_4 + 36 u_2^2
t_4^2 + 54 u_2^2 u_3 t_5 - 80 u_2^3 t_6) t_6\\ &\qquad + (8 u_2 u_3^3 -
30 u_2^2 u_3 t_4 + 50 u_2^3 t_5) t_7,\\
\intertext{of total degree 18, and}
G &= u_3^2 t_4^2 t_5 - 4 u_2 t_4^3 t_5 - 4 u_3^3 t_5^2 + 18 u_2 u_3 t_4
t_5^2 - 27 u_2^2 t_5^3\\ &\qquad + (2 u_3^3 t_4 - 8 u_2 u_3 t_4^2 - 6 u_2
u_3^2 t_5 + 36 u_2^2 t_4 t_5 - 8 u_2^2 u_3 t_6) t_6\\ &\qquad + (16
u_3^4 - 74 u_2 u_3^2 t_4 + 40 u_2^2 t_4^2 + 110 u_2^2 u_3 t_5 - 200
u_2^3 t_6) t_7,
\end{align*}
of total degree 19.  Let $P \mapsto \tild P$ denote the operation of
expanding the $t_i$ in terms of $u_2$, $u_3$, and $\sigma_1$, \dots,
$\sigma_n$ as specified by Corollary~\ref{c:cyclic}.  This operation
preserves total degree.  Both $\tild F$ and $\tild G$ have degree $6$
in~$u_3$.  Their resultant with respect to $u_3$ therefore has total
degree $6\times 19 = 114$, and it must have the polynomial $\alpha_7$ of total
degree $2\Delta_7 = 76$ as a factor.

The resultant $\Res(\tild F, \tild G)$ seems to be too large to compute
and factor explicity, but we can describe its unwanted factors as
follows with a little computer assistance.  First observe that every
term in $F$ and $G$ is divisible by either $u_2$ or $u_3$, and hence the
same is true of $\tild F$ and $\tild G$.  It follows that $\Res(\tild F,
\tild G)$ is divisible by $u_2$.  In fact $u_2^7 \divides \Res(\tild F,
\tild G)$, as we will see in Lemma~\ref{l:mplcty} below. Next, consider
the polynomials
\begin{align*}
F_1 &= 4 u_3^3 - 15 u_2 u_3 t_4 + 25 u_2^2 t_5,\\
G_1 &= 7 u_3^2 t_4 - 20 u_2 t_4^2 - 5 u_2 u_3 t_5 + 100 u_2^2 t_6,
\end{align*}
which are closely related to the coefficients of $t_7$ in $F$ and $G$.
Specifically, $F_1 = (2u_2)^{-1} \partial F/\partial t_7$ and $G_1 =
u_2^{-1}(2 u_3 F_1 - \partial G/\partial t_7)$.  We will show that
$\Res(\tild{F}_1, \tild{G}_1)$ divides $\Res(\tild F, \tild G)$.

First we claim that if $F_1 = G_1 = 0$, then $F = G = 0$.  The ideal
$\langle F_1,G_1 \rangle$ does not contain $F$ and $G$, but by some easy
calculations, it does contain $u_2 F$, $u_3 F$, $u_2 G$, and $u_3 G$.
If $F_1=G_1=0$, then all four of these polynomials vanish; so if either
$u_2\ne 0$ or $u_3\ne 0$, we must have $F=G=0$, while if $u_2 = u_3 =
0$, we already know that $F=G=0$.  This establishes the claim.  It
follows that $\tild{F}_1 = \tild{G}_1 = 0$ implies $\tild F = \tild G =
0$.  Consequently, wherever $\Res(\tild{F}_1,\tild{G}_1)$ vanishes, so
does $\Res(\tild F,\tild G)$.  Algebraically, this means that every
irreducible factor of $\Res(\tild{F}_1,\tild{G}_1)$ divides $\Res(\tild
F,\tild G)$.  The resultant of $\tild{F}_1$ and $\tild{G}_1$ with
respect to~$u_3$ is simple enough to compute explicitly.  It has total
degree 30, and it factors as $u_2^3$ times an irreducible polynomial in
$\Qfield[u_2,\sigma_1,\ldots,\sigma_7]$ of total degree~24.

Thus, not only does $\Res(\tild{F}_1,\tild{G}_1)$ divide $\Res(\tild
F,\tild G)$, but $u_2^4 \Res(\tild{F}_1,\tild{G}_1)$ does also.  The
quotient by the latter polynomial has total degree $114 - 8 - 30 = 76 =
2\Delta_7$, so it must be a scalar multiple of the desired polynomial
$\alpha_7$, $\beta_8$, or $\beta^*_8$; there are no more unwanted
factors.  The scalar can be computed by setting $\sigma_2$, \dots,
$\sigma_7$ to zero, and we find that
\begin{equation*}
\frac {2^{101} 5^5 \Res(\tild F,\tild G,u_3)} {u_2^4
\Res(\tild{F}_1,\tild{G}_1,u_3)}
\end{equation*}
is $\alpha_7$, $\beta_8$, or $\beta^*_8$, according to whether the
crossing parity $\epsilon$ is $0$, $+1$, or~$-1$.  It remains only to
prove the following lemma.

\begin{lemma}\label{l:mplcty}\sl
With the definitions above, $u_2^7 \divides \Res(\tild F, \tild G, u_3)$.
\end{lemma}

\begin{proof}[Sketch of proof]
By direct calculation on a computer, $u_2^7$ divides $\Res(F,G,u_3)$
but $u_2^8$ does not.  The only component of $V(F,G)$ lying on the
hyperplane $u_2=0$ is the linear variety $V(u_2,u_3)$, so $V(F)$ and
$V(G)$ must intersect with multiplicity~7 along $V(u_2,u_3) \subset
\Aspace^6$.  Now, assuming $n=7$ for definiteness, pull back via the
projection $\pi: \Aspace^9 \to \Aspace^6$ that maps $(u_2,u_3, \sigma_1,
\ldots, \sigma_7) \mapsto (u_2,u_3, t_4,\ldots,t_7)$.  Because $\pi$ is
smooth, the intersection multiplicity of $V(\tild F)$ and $V(\tild G)$
along $\pi^{-1}V(u_2,u_3) = V(u_2,u_3) \subset \Aspace^9$ is also~$7$.
For fixed generic values of $\sigma_1$, \dots, $\sigma_7$, we therefore
have $u_2^7 \divides \Res(\tild F,\tild G,u_3)$.  We conclude that
this divisibility holds globally as well.
\end{proof}

For the rest of this section, we turn our attention to semicyclic
$(n+1)$-gons with $n=3$, $4$, $5$, and~$6$.  To state the area formulas
most cleanly we introduce a notion of parity for semicyclic polygons.
Let $n$ be even, and observe that the quantities $e_1 = \tau_1 +
\tau_{n-1} = \sum q_i - \sum q_i^{-1}$ and $\half e_{n/2} = \sum
(-1)^i \tau_i = \prod (1-q_i)$ are both pure imaginary.  (Compute their
complex conjugates using $\bar{q}_i = q_i^{-1}$.) Hence their product is
real.  Let $\epsilon \in \{-1,0,+1\}$ be its sign.  Then we have
\begin{align*}
\epsilon \abs{K} \sqrt{\sigma_n} &= \epsilon \cdot \qtr r^2 \abs{\tau_1
+ \tau_{n-1}} \cdot \abs{v_1 - v_2} \abs{v_2 - v_3} \cdots \abs{v_n -
v_1}\\ &= \epsilon \cdot \qtr r^2 \abs{\tau_1 + \tau_{n-1}} \cdot r^n
\abs{1-q_1} \cdots \abs{1-q_n}\\ &= r^{n+2} \cdot \qtr e_1 \cdot \half
e_{n/2}.
\end{align*}
Define $w = 2\epsilon \abs{K} \sqrt{\sigma_n} = \epsilon \sqrt{u_2 t_n}$
for $n$ even, and let $w = 0$ for $n$ odd.  Our formulas for $\alpha'_4$
and $\alpha'_6$ will factor when written in terms of $w$ rather than
$\sigma_n$.  We do not know whether this type of factorization occurs in
general.

For $n\in \{3,4\}$, it is simplest to use the main identity directly.
Defining $e_2=0$ if $n=3$, we have
\begin{align*}
\qtr E(r^2 x)^2 &= \qtr r^4 e_1^2 x^2 + \half r^6 e_1 e_2 x^3 + \qtr r^8
e_2^2 x^4\\ &= u_2 x^2 + 2w x^3 + \cdots
\end{align*}
so, by Theorem~\ref{t:main}, the cubic $1 - \sigma_1 x + (\sigma_2 +
u_2) x^2 - (\sigma_3 - 2w) x^3$ factors as $-(r^2 x - \qtr) D(r^2
x)^2$.  In particular, its discriminant vanishes.  Replacing $x$
by~$-x^{-1}$, we recover equation (\ref{e:semiquad}) for $n=3$, and for
$n=4$ we have factored $\alpha'_4$ as the product of two discriminants
$\beta'_4$ and $(\beta'_4)^*$ corresponding to $\eps=+1$ and $\eps=-1$
respectively.

For larger $n$ we need Corollary~\ref{c:semi}.  For $n=5$, it says that
$u_2 + t_3 z + t_4 z^2 + t_5 z^3$ is the square of the linear polynomial
$E(r^2 z) / (2z)$, which yields the two equations $t_3^2 - 4 u_2 t_4 =
0$ and $t_5 = 0$.  Their degrees in $r^2$ are $6$ and $5$ respectively,
so their resultant with respect to $r^2$ has the correct total degree
$2\Delta'_5 = 30$.  (Remember that $r^2$ has degree~$1$.)  It remains
only to scale the resultant to be monic in $-4u_2$, and we get
$$\alpha'_5 = \qtr \Res(t_3^2 - 4 u_2 t_4, t_5, r^2).$$
For $n=6$, Corollary~\ref{c:semi} gives us the factorization
$$u_2 + t_3 z + t_4 z^2 + t_5 z^3 + t_6 z^4 = \qtr r^4
(e_1 + e_2 r^2 z + e_3 r^4 z^2)^2.$$
Using $w = \qtr r^8 e_1 e_3$, we derive the equations
\begin{align*}
u_2 t_5 - t_3 w &= 0,\\
u_2 + t_3 z + (t_4 - 2w) z^2 &= \qtr r^4 (e_1 + e_2 r^2 z)^2,
\end{align*}
the second of which implies that $t_3^2 - 4 u_2 (t_4 - 2w) = 0$.  Thus
we can form the resultant of $t_3^2 - 4 u_2 (t_4 - 2w)$ and $u_2 t_5 -
t_3 w$ to eliminate $r^2$ and obtain a multiple of the desired area
formula.  The resultant is small enough to compute and factor
symbolically, and we obtain $\alpha'_6 = (\beta'_6)(\beta'_6)^*$ where
\begin{equation*}
\beta'_6 = \frac{\Res(t_3^2 - 4 u_2 (t_4 - 2\sqrt{u_2 t_6}), u_2 t_5 -
  t_3 \sqrt{u_2 t_6}, r^2)} {4 u_2^6},
\end{equation*}
and $(\beta'_6)^*$ is $\beta'_6$ with the opposite sign on
$\sqrt{u_2t_6}$.

\section{Degree Calculations}\label{s:radius}

In this section we show by elementary means that the homogeneous polynomials
$\alpha_n$ and $\alpha'_n$ have degrees $2\Delta_n$ and $2\Delta'_n$ respectively.
First we explain why the degrees cannot be smaller.
In \cite{robbins}, Robbins shows that $\deg(\alpha_n) \ge 2\Delta_n$ by constructing $\Delta_n$ cyclic $n$-gons with generically different squared areas from a given set
of edge lengths.  He takes the edge lengths to be nearly equal if $n$
is odd, and adds a much shorter edge if $n$ is even.  For semicyclic polygons,
we can take the edge lengths to be nearly equal if $n$ is even;
the argument of \cite{robbins} then yields the desired number $\Delta'_n$ of
semicyclic $n$-gons.

Suppose now that $n$ is odd.
It is not necessary (and in fact not possible) to construct $\Delta'_n$
inequivalent semicyclic polygons with given positive real edge lengths $a_j$.
It suffices instead to construct $\Delta_n$ configurations
$(r, q_1, \ldots, q_n)$ of complex numbers satisfying
\begin{equation}
\label{e:qfromr} a_j^2 = r^2 (2 - q_j - q_j^{-1}), \qquad j=1,\ldots,n;
\end{equation}
and $q_1 \cdots q_n = -1$, since it is from these equations,
together with the relation $16K^2 =
-r^4 (\sum q_j - \sum q_j^{-1})^2$, that one derives the existence and uniqueness of
the irreducible polynomial $\alpha'_n$.  In our configurations $r$ is always
real and positive, but sometimes $r < \min \{a_j/2\}$, in which case the $q_j$
are negative real numbers instead of complex numbers of norm~$1$.  The
plan is to regard each $q_j$ as a function of $r$ by choosing a branch of
equation (\ref{e:qfromr}), and then find values of $r$ such that $q_1
\cdots q_n = -1$.  

Let $n=2m+1$, let the first $2m$ edge lengths be large and nearly equal,
and let $a_n = 2$.  To find configurations with
$r > \max \{a_j/2\}$, choose arbitrarily whether $0 < \arg q_n < \pi$
(the short edge goes ``forward'') or $-\pi < \arg q_n < 0$ (``backward''),
and likewise choose a set of $k < m$ of the long edges to go backward.  
Then there exist $m-k$ semicyclic polygons with the given edge lengths and
edge directions whose angle sums
$\sum \arg q_j$ are $\pi$, $3\pi$, \dots, $(2m-2k-1)\pi$.  (Apply the
Intermediate Value Theorem to $\sum \arg q_j$ as $r$ varies from
$\max \{a_j/2\}$ to~$\infty$.)  The total number of such configurations is
$$\sum_{k=0}^{m-1} 2\binom{2m}{k}(m-k) = m\binom{2m}{m}.$$

To find configurations
with $r < \min \{a_j/2\} = 1$, choose the branch $q_j < -1$ for exactly $m$ of
the long edges, and choose the branch $q_j > -1$ for the other long edges.
Let $\eps_j = +1$ or $\eps_j = -1$ respectively.
As $r\to 0$, the product $q_1 \cdots q_{2m}$ approaches the
constant $\prod_{j=1}^{2m} a_j^{2\eps_j}$, and hence $q_1 \cdots q_n$
approaches $0$ if $q_n > -1$ or $-\infty$ if $q_n < -1$.  By choosing the branch
for $q_n$ according
to whether $q_1 \cdots q_{2m}$, evaluated at $r=1$, is greater or less than $1$,
we guarantee that $q_1 \cdots q_n = -1$ for some intermediate value of~$r$.
Thus we obtain another $\half \binom{2m}{m}$ configurations.  (The factor of $1/2$
is present because inverting every $q_j$ preserves the radius and the squared area;
it corresponds to reversing the orientation.)
The total number of configurations is therefore at least
$$(m + \half)\binom{2m}{m} = \frac{n}{2} \binom{n-1}{\floor{\frac{n-1}{2}}} = \Delta'_n.$$

To establish matching upper bounds on the degrees of $\alpha_n$ and $\alpha'_n$, we
proceed indirectly.  First we revive an argument of M\"obius from the 19th century \cite{moebius}, which produces a polynomial of degree $\Delta_n$ that relates $r^2$
for a cyclic polygon to the squared side lengths.  Hence there are generically at most
$\Delta_n$ circumradii for a given set of edge lengths.  For generic side lengths (in particular, no two equal) and a radius $r$ that admits a solution $(q_1,\ldots,q_n)$
to the system of equations (\ref{e:qfromr}) and $q_1\cdots q_n = 1$,
the solution is unique up to inverting all the $q_j$.  (Any other solution would
differ by inverting a proper subset of the $q_j$, so those $q_j$ would need to have
product $\pm 1$.)  Thus, because $r$ and the $q_j$ determine the area, there are generically
at most $2\Delta_n$ possible signed areas, so $\deg(\alpha_n) \le 2\Delta_n$.
The same argument applied to semicyclic polygons will yield
$\deg(\alpha'_n) \le 2\Delta'_n$.

Given a cyclic polygon with circumradius $r$ and side lengths $2y_1$, \dots, $2y_n$,
let $\theta_j = \sin^{-1}(y_j/r)$ be half the angle subtended by the $j\th$ side.
Then, for some choice of signs $\eps_2$, \dots $\eps_n$ (namely, $\eps_j$ is
$+1$ or $-1$ according to whether the $j\th$ side goes ``forward'' or ``backward''
relative to the first side), the sum $\theta_1 + \epsilon_2 \theta_2 + \cdots +
\epsilon_n \theta_n$ is a multiple of $\pi$.  Therefore
\begin{equation}\label{e:moebius}
\prod_{\eps_j = \pm 1} r^n \sin(\theta_1 + \epsilon_2 \theta_2 + \cdots +
\epsilon_n \theta_n) = 0.
\end{equation}
The factors of $r$ make this a polynomial relation over $\Qfield$ between $r^2$
and the squared side lengths.  To see why,
introduce the variables $x_j = r\cos \theta_j = (r^2 - y_j^2)^{1/2}$ and rewrite
(\ref{e:moebius}) as
\begin{equation}\label{e:mexpand}
\prod_{\eps_2, \ldots, \eps_n} \frac{1}{2i} \left[ \prod_{j=1}^n (x_j
+ i \eps_j y_j) - \prod_{j=1}^n (x_j - i \eps_j y_j) \right] = 0
\end{equation}
where $\eps_1 = 1$.  The left-hand side of (\ref{e:mexpand}) has a great deal of
symmetry.  Obviously, flipping the sign of $y_j$ is equivalent to negating $\eps_j$.
Flipping the sign of any $x_j$ is equivalent to flipping $\eps_j$ and negating each
product over~$j$.  If $j=1$, we can restore the condition $\eps_1=1$ by flipping every 
$\eps_j$ and negating every bracket.  All these operations just permute and possibly
negate all the $2^{n-1}$ bracketed factors, so they leave the overall expression unchanged.
Therefore, in the expansion of (\ref{e:mexpand}), each $x_j$ occurs only to even powers,
and hence each $x_j^2$ can be replaced by $r^2 - y_j^2$.  Likewise each $y_j$ occurs
only to even powers.  Thus we obtain a polynomial equation $M(r^2,y_1^2,\ldots,y_n^2)=0$.

The remaining part of M\"obius' argument uses series expansion to find the
degrees of the leading and trailing terms of $M$.  Fix
the $\eps_j$, and rewrite the bracketed factor of (\ref{e:mexpand}) as
$$\prod_{j=1}^n \Bigl( \sqrt{r^2 - y_j^2} + i \eps_j y_j \Bigr)
- \prod_{j=1}^n \Bigl( \sqrt{r^2 - y_j^2} - i \eps_j y_j \Bigr).$$
To find the term of highest degree in $r$, expand around $r=\infty$; the highest terms cancel,
so the degree is $n-1$.  To find the term of lowest degree, expand around $r = 0$ to get
$$\prod_{j=1}^n i y_j \left( 1+\eps_j  - \frac{r^2}{2y_j^2} - \cdots \right)
- \prod_{j=1}^n i y_j \left( 1-\eps_j  - \frac{r^2}{2y_j^2} - \cdots \right).$$
Its initial term has degree $\min(k,n-k)$ in $r^2$, where $k$ is the number of
$\eps_j$ equal to $-1$.  Therefore $M$ is a power of $r^2$ times a polynomial in
$r^2$ of degree
$$2^{n-1} \frac{n-1}{2} - \sum_{k=0}^{n-1} \binom{n-1}{k} \min(k,n-k),$$
which simplifies to $\Delta_n$.  We can factor out the unwanted power of $r^2$
because it was not needed to make equation (\ref{e:moebius}) hold.

For semicyclic polygons, the signed sum of the $\theta_j$ is an odd multiple
of $\pi/2$.  Equation (\ref{e:moebius}) therefore becomes
$$\prod_{\eps_j = \pm 1} r^n \cos(\theta_1 + \epsilon_2 \theta_2 + \cdots +
\epsilon_n \theta_n) = 0,$$ which expands to a polynomial relation
$M'(r^2, y_1^2, \ldots, y_n^2) = 0$.  Using series expansion again, one finds
that $M'$ is monic of degree $n2^{n-2}$ in $r^2$, and its lowest nonzero term
has the same degree as that of $M$.  Hence $M'$ is a power of $r^2$ times a
polynomial whose degree in $r^2$ is $\Delta_n + 2^{n-2} = \Delta'_n$.

\section{Specializations}\label{s:special}

Corollaries \ref{c:cyclic} and~\ref{c:semi}, which relate the generalized
Heron polynomials $\alpha_n$ and $\alpha'_n$ to the polynomials $P_n(z)$ and
$P'_n(z)$, allow us to understand and factor certain specializations of $\alpha_n$ and
$\alpha'_n$.  With a little extra work one can describe some of the factors
explicitly.  In this section we offer two such results concerning
cyclic $n$-gons with $n$ odd.


Let $n=2m+1 \ge 5$, and consider the constant term of $\alpha_n$; that is, let
$u_2=0$.  Then $P_n(z)$ has $m-1$ double roots if and only if
either $(P_n|_{u_2=0})/z$ has $m-1$ double roots, or $u_3=0$ and
$(P_n|_{u_2=u_3=0})/z^2$ has $m-2$ double roots.  Geometrically, the projective
variety
\begin{multline*}
X = \bigl\{\, [u_2:\cdots:u_m:t_{m+1}:\cdots:t_{2m+1}] \bigm| \\
\text{$P_n(z)$ factors as $(b_0+b_1z)(c_0+c_1z+\cdots+c_{m-1}z^{m-1})^2$} \,\bigr\}
\end{multline*}
intersects the hyperplane $\{u_2=0\}$ in two irreducible components, one
corresponding to $b_0=0$ and one corresponding to $c_0=0$.
The second component has intersection multiplicity
two because $X$ is tangent to $\{u_2=0\}$ along it.  Chasing through the geometric
interpretation of $\alpha_n$ (after Corollary~\ref{c:cyclic}), we find that
$\alpha_n|_{u_2=0}$, as a polynomial in $\sigma_1$, \dots,~$\sigma_n$, is an irreducible polynomial times the square of another irreducible.  The factors are not necessarily irreducible as polynomials in the side lengths $a_i$, however.

\begin{proposition}\label{p:const}\sl If $n$ is odd, the constant term
of $\alpha_n$ factors as
$$\alpha_n|_{16K^2=0} = \gamma_n^2 \prod (a_1 \pm a_2 \pm \cdots \pm a_n)$$
where the product is over all $2^{n-1}$ sign patterns.
\end{proposition}

\begin{proof}
Heron's formula takes care of the case $n=3$, so we may assume $n\ge 5$ and apply
the analysis above.  By Corollary~\ref{c:cyclic}, cyclic $n$-gons satisfy
$$P_n(z) = (\qtr-r^2 z) \bigl( D(r^2z)/z \bigr)^2,$$
so the factor $\gamma_n^2$ corresponds to $d_1 = 0$,
and the other factor corresponds to $[\qtr:-r^2] = [0:1]$ and represents projective
solutions at $r^2=\infty$.  The presence of the linear factors $a_1 \pm a_2 \pm \cdots
\pm a_n$ in the constant term was proved in~\cite{varfolomeev}, and they
correspond to solutions with $r^2=\infty$: As a signed sum of edge lengths approaches
zero, the polygon can degenerate to a chain of collinear line segments,
which has zero area and infinite circumradius.  (One can easily construct
a curve of solutions to the equations (\ref{e:vq}) tending to any such
point at infinity.)  For $n\ge 3$, the product of these $2^{n-1}$ linear factors is
symmetric in the $a_i^2$, so by irreducibility, no other factors can appear.
\end{proof}

The same kind of analysis applies to $\alpha_n$ when the side length $a_n$
goes to zero, and so $t_n = \sigma_n = 0$.  It's geometrically clear that
the result should be divisible
by $\alpha_{n-1}^2$ (as the $n\th$ side shrinks to zero,
it can go either ``forward'' or ``backward''), and the algebra confirms it.
The intersection of $X$ with the
hyperplane $t_n = 0$ includes a component of multiplicity two where 
$$t_{n-1} = -\sigma_{2m} + \qtr t_m^2 = 0$$ and $P_n(z)$, considered as
degree~$2m-3$, has $m-2$ double roots.  Substituting the solutions $t_m = \pm 2 \sqrt{\sigma_{2m}}$ back into the definitions of $t_{m+1}$ through $t_{2m-1}$, 
we recover the definitions of the $t_j$ for $n=2m$ and $\eps=\pm 1$
(see Corollary~\ref{c:cyclic}), and observe that $u_m$ becomes $t_m$.
Thus $P_n(z)$ specializes to $P_{n-1}(z)$, and so $\alpha_n |_{a_n=0}$ is
divisible by $(\beta_{n-1})^2 (\beta^*_{n-1})^2 = \alpha_{n-1}^2$.

The other component of $X \cap \{t_n=0\}$ corresponds to solutions in which the
the leading coefficient $r^2$ of the linear factor $(\qtr - r^2 z)$
vanishes.  We can describe these solutions
explicitly.

\begin{proposition}\label{p:degen}\sl If $n\ge 5$ and $n$ is odd, then
$$\alpha_n|_{a_n=0} =  \alpha_{n-1}^2 \prod 
\bigl( 16K^2 + (a_1^2 \pm a_2^2 \pm \cdots \pm a_{n-1}^2)^2 \bigr),$$
the product taken over all sign patterns with $(n-1)/2$ minus signs.
\end{proposition}


\begin{proof}
It suffices to show that all the factors in the product are present;
the result then follows by comparing degrees (both sides being monic in $16K^2$).

Fix generic positive real values for $a_1$, \dots, $a_{n-1}$ and signs $\eps_1$,
\dots,~$\eps_{n-1} \in \{-1,+1\}$ with $\eps_1 = +1$ and $\sum \eps_j = 0$.
Recall that cyclic polygons satisfy
$$q_j^2 + (a_j^2/r^2 - 2) q_j + 1 = 0, \qquad 1\le j\le n,$$
and $q_1 q_2\cdots q_n = 1$, and any solution to
these $n+1$ equations also satifies $\alpha_n$ with the squared area given by
equation (\ref{e:area}).  For sufficiently small positive values of~$a_n$, the
method of $\S$5 shows that there
exist solutions $(r^2,q_1,\ldots,q_n)$ with $r^2 < 0$ and
$q_j \approx (a_j^2/r^2)^{\eps_j}$ for $1\le j < n$; furthermore
$a_n^2/r^2$ tends to a constant.  Equation~(\ref{e:area}) now implies
\begin{equation*}
\lim_{a_n\to 0} 16K^2 =
-\biggl(\sum_{\eps_j=+1} a_j^2 - \sum_{\eps_j=-1} a_j^2 \biggr)^2,
\end{equation*}
so the set of solutions with $a_n^2 = r^2 = 0$ includes all points on the hypersurface
$16K^2 + \bigl(\sum_{j=1}^{n-1} \eps_j a_j^2\bigr)^2 = 0$.
\end{proof}

\section{Conclusions}

The generalized Heron polynomials $\alpha_n$ for 
cyclic $n$-gons and $\alpha'_n$ for semicyclic $(n+1)$-gons are
defined implicitly by $n+2$ equations in $n+1$ unknowns:
$n$~equations (\ref{e:vq}) relating the side lengths to the
vertex quotients $q_1$, \dots,~$q_n$ and the radius $r$;
equation (\ref{e:area}) which expresses
the squared area $K^2$, or equivalently $u_2 = -4K^2$, in the same way;
and the equation $q_1 q_2 \cdots q_n = \delta$.
Our analysis eliminates the variables $q_j$ (and in the cyclic case,
$r$ also) at the cost of introducing
$\floor{(n-5)/2}$ unwanted quantities $u_3$, \dots, $u_m$.
The reduction in the number of auxiliary variables allows us,
for small $n$, to eliminate them by ad hoc means and obtain
formulas for $\alpha_n$ and $\alpha'_n$.

The quantities $u_2$, \dots, $u_m$ appear on equal footing in this
analysis, so we could equally well eliminate all but $u_k$ for some
$k>2$ and obtain a polynomial relation, presumably of total degree
$k\Delta_n$ or $k\Delta'_n$, between $u_k$ and the squares of the side lengths.
Unfortunately, we do not yet have a geometric interpretation for $u_3$
or the higher $u_k$.

For large $n$ the goal of eliminating $u_3$, \dots, $u_m$ seems rather
distant, but Corollaries \ref{c:cyclic} and~\ref{c:semi} still illuminate
aspects of the polynomials $\alpha_n$ and $\alpha'_n$.
In particular, Corollary \ref{c:cyclic} establishes a close relationship
between $\alpha_{2m+1}$ and $\alpha_{2m+2}$, generalizing those between
Heron's and Brahmagupta's formulas and between Robbins' pentagon and
hexagon formulas.

It may be of some interest to know how our main results were obtained.
Robbins solved many combinatorial and algebraic problems in his lifetime
by what he called the ``Euler method'': calculate examples, using a
computer if convenient; discover a general pattern; and prove it, with hints from
further calculations if necessary.  His work on generalized Heron polynomials
took this approach but was somewhat frustrated by lack of data from which
to generalize.  As explained in \cite{robbins}, Robbins first found $\alpha_5$ and
the closely related polynomial $\alpha_6$ by interpolating from several
dozen numerical examples, and he rewrote them concisely in terms of variables $t_i$
(slightly different from ours) by interacting with a computer algebra system.
Neither step is particularly feasible for $\alpha_7$, whose expansion in terms
of the symmetric functions $\sigma_k$ has almost a million coefficients.
It was possible, however, to evaluate certain specializations of $\alpha_7$ and
$\alpha_8$ by interpolation, and to conjecture Propositions \ref{p:const}
and~\ref{p:degen}.  For instance, we discovered that the constant term of
$\alpha_7$ (the specialization $u_2=0$) is divisible by the square of the
discriminant of $$t_4 + t_5 z + t_6z^2 + t_7 z^3,$$ where the $t_i$ are defined
as in Corollary~\ref{c:cyclic} but with $u_2=u_3=0$.  Unfortunately, the
hidden presence of $u_3$ made it difficult to guess the rest of $\alpha_7$.

We introduced semicyclic polygons and their area polynomials $\alpha'_n$ in an effort
to obtain more data to study.  In particular, we noticed that
the mysterious cubic discriminant so prominent in $\alpha_5$ appeared already in
the simpler polynomial $\alpha'_3$, and we hoped
that whatever new phenomenon arose in $\alpha_7$ would also appear in 
$\alpha'_5$, which we could compute by interpolation.  In fact $\alpha_7$ turns out to be
rather different from $\alpha'_5$, but it was the struggle to simplify $\alpha'_5$
that led us to manipulate the relations between the $\sigma_i$ and $\tau_j$ by hand
and thence to discover the main identity (Theorem~\ref{t:main}).  In the end, the
most crucial calculations turned out to be those we did on the blackboard.

\section*{Postscript}

David Robbins (1942--2003) had an exceptional ability to see and communicate
the simple essence of complicated mathematical issues, and to discover elegant
new results about seemingly well-understood problems.  He taught and 
inspired a long sequence of younger mathematicians including the two surviving
authors.  His interest in cyclic polygons began at age 13 when he derived a
version of Heron's formula.  In the early 1990's he discovered the area
formulas for cyclic pentagons and hexagons.  When diagnosed with a terminal
illness in the spring of 2003, he chose to work on this topic once again. 
Sadly, he did not live to see the discovery of the main identity or the
heptagon formula.  This paper is dedicated to the memory of our friend and
colleague, whose loss is keenly felt.

\section*{Acknowledgements}

We thank David Lieberman and Lynne Butler for helpful discussions, and Joel
Rosenberg for simplifying the proof of the main identity.

\bibliography{heron}
\bibliographystyle{plain}

\end{document}